\documentclass[12pt, twoside]{article}
\usepackage{mathrsfs, amsmath,amssymb,amsfonts, amsthm, dsfont}
\usepackage{times}
\usepackage{enumerate}

\def\titlerunning#1{\gdef\titrun{#1}}
\makeatletter
\def\author#1{\gdef\autrun{\def\and{\unskip, }#1}\gdef\@author{#1}}
\def\address#1{{\def\and{\\\hspace*{18pt}}\renewcommand{\thefootnote}{}%
\footnote {#1}}%
\markboth{\autrun}{\titrun}}
\makeatother
\def\email#1{e-mail: #1}
\def\subjclass#1{{\renewcommand{\thefootnote}{}%
\footnote{\emph{Mathematics Subject Classification (2010):} #1}}}

\newtheorem{thm}{Theorem}[section]

\newtheorem{lem}[thm]{Lemma}

\theoremstyle{definition}

\numberwithin{equation}{section}

\begin{document}


\baselineskip=17pt


\titlerunning{Sharp bounds for the intersection of nodal lines with certain curves}

\title{Sharp bounds for the intersection of nodal lines with certain curves}

\author{Junehyuk Jung}

\date{\today}

\maketitle

\address{J. Jung: Princeton University, Princeton, NJ 08544; \email{junehyuk@math.princeton.edu}}

\subjclass{Primary XXXX; Secondary YYYY}


\begin{abstract}
Let $Y$ be a hyperbolic surface and let $\phi$ be a Laplacian eigenfunction having eigenvalue $-1/4-\tau^2$ with $\tau>0$. Let $N(\phi)$ be the set of nodal lines of $\phi$.  For a fixed analytic curve $\gamma$ of finite length, we study the number of intersections between $N(\phi)$ and $\gamma$ in terms of $\tau$. When $Y$ is compact and $\gamma$ a geodesic circle, or when $Y$ has finite volume and $\gamma$ is a closed horocycle, we prove that $\gamma$ is ``good'' in the sense of \cite{pa}. As a result, we obtain that the number of intersections between $N(\phi)$ and $\gamma$ is $O(\tau)$. This bound is sharp.
\end{abstract}

\section{Introduction}
Let $Y$ be a hyperbolic surface (a Riemannian surface with constant curvature $-1$) and let $\triangle$ be the Laplace-Beltrami operator on $Y$. Let $\phi \in L^2(Y)$ be a real valued eigenfunction for $\triangle$ with eigenvalue $-(1/4+\tau^2)$, that is
\[
\triangle \phi =-(\frac{1}{4}+\tau^2) \phi.
\]
Let $N(\phi)$ denote the zero set of $\phi$ on $Y$. Since $Y$ is an analytic manifold and $\phi$ is an analytic function, $N(\phi)$ is a union of analytic curves on $Y$. Therefore, for an analytic segment $\gamma$ of finite length, either $\phi$ vanishes on $\gamma$ or $\gamma \cap N(\phi)$ consists of a finite number of points.

It has been shown in several ways that $N(\phi)$ shares common properties with the zero set of a polynomial of degree $\tau$ for large $\tau>0$ (see \cite{df}, \cite{sa}, and \cite{zeld}). Accordingly we expect that
\[
|\gamma \cap N(\phi)| \ll_{Y,\gamma} \tau 
\]
will hold for all $\phi$ which does not vanish on $\gamma$. Here and elsewhere $A \ll_\omega B$ means $|A| < CB$ for some constant $C>0$ depending only on $\omega$. Unless otherwise stated, in what follows, we assume that $\tau>0$ is sufficiently large (say $\tau >100$.) 

The main result of the article is:
\begin{thm}\label{theorem}
Assume either one of the following:
\begin{itemize}
\item[-] $Y$ is a compact surface and $\gamma$ is a geodesic circle,
\item[-] $Y$ is a noncompact surface of finite volume and $\gamma$ is a closed horocycle.
\end{itemize}
Then $N(\phi)\cap \gamma$ consists of a finite number of points and we have
\begin{equation}\label{equation123}
|N(\phi) \cap \gamma| \ll_{Y,\gamma} \tau.
\end{equation}
\end{thm}
We remark that, in spite of its appearance, (\ref{equation123}) is not a local feature of the eigenfunctions on hyperbolic surfaces. That is, without any assumptions on $Y$, the assertions in Theorem \ref{theorem} are not true in general. We give two examples to emphasize the role of the global geometry of $Y$.

\begin{enumerate}
\item Let $Y=\Gamma_0\backslash \mathbb{H}$ where $\Gamma_0 = \{\left( \begin{array}{cc} 1 & n \\ 0 & 1 \end{array} \right)| n \in \mathbb{Z}\}$. Let the horocycle $\gamma$ be given by $\{x+i| x\in \mathbb{R}\}$. Consider a one-parameter family of eigenfunctions
\[
f_\tau(x+iy)= \sqrt{y}K_{i\tau}(2\pi y)\cos 2\pi  x.
\]
Since there are infinitely many $\tau \in \mathbb{R}$ such that $K_{i\tau}(2\pi)=0$, we can find a sequence of eigenfunctions $f_{\tau_n}$ with $\tau_n \to \infty$ such that they all vanish on $\gamma$.

\item Let $Y$ be as above and let $\gamma$ be an analytic curve which is not part of a vertical line. For $n>0$, let
\[
f_n(x+iy) = \sqrt{y}K_i(2\pi n y) \cos 2\pi  n x
\]
be a sequence of eigenfunctions having eigenvalue $-5/4$. Then for all sufficiently large $n$, $N(f_n)$ does not vanish identically on $\gamma$ and we have
\[
|N(f_n) \cap \gamma| \gg n.
\]
\end{enumerate}

In fact, if one assumes $Y$ has finite volume, then a nonzero eigenfunction $\phi$ in $L^2(Y)$ cannot vanish either on a horocycle or on a geodesic circle. Let us sketch a proof for the horocycle case; the other case can be treated similarly.

Without loss of generality, assume that $Y=\Gamma\backslash \mathbb{H}$ where $\Gamma$ is reduced at $\infty$. We further assume that $\gamma$ is given by $\{x+iy_0| x\in \mathbb{R}\}$ for some $y_0>0$. Say the Fourier expansion of the eigenfunction $\phi \in L^2(Y)$ at $\infty$ is given by
\[
\phi(x+iy)=\sum_{n \neq 0} a_{\phi}(n) \sqrt{y}K_{i\tau}(2\pi |n| y )e^{2\pi i n x}.
\]
If $\phi$ vanishes on $\gamma$, then
\[
\int_0^1 |\phi(x+iy_0)|^2 dx = \sum_{n \neq 0} |a_{\phi}(n)|^2 y_0 K_{i\tau}(2\pi |n| y_0 )^2 =0.
\]
Since $K_{i\tau}(2\pi |n| y_0 )=0$ for only finitely many $n$, this implies that $a_\phi(n)=0$ for all but finitely many $n$. Now recall that the horocycle becomes equidistributed on $Y$ when one increases the length of the horocycle \cite{horo}. (For the geodesic circle, the first equidistribution result is due to Delsarte \cite{geod}. We give a short proof of the result when $Y$ is compact in the appendix.) Therefore, we have
\begin{align*}
||\phi||_{L^2(Y)}^2 &= \lim_{y \to 0}\int_0^1 |\phi(x+iy)|^2 dx\\
 &= \lim_{y \to 0} \sum_{n \neq 0} |a_{\phi}(n)|^2 y K_{i\tau}(2\pi |n| y )^2.
\end{align*}
By interchanging the limit and the finite summation, we obtain
\begin{align*}
 &=\sum_{n \neq 0} |a_{\phi}(n)|^2 \lim_{y \to 0}   y K_{i\tau}(2\pi |n| y )^2\\
 &=0
\end{align*}
which contradicts the assumption $\phi \neq 0$.

It is noteworthy that Theorem \ref{theorem} does not make any arithmetic assumptions. With such an assumption, when $\gamma$ is a closed horocycle for instance, the following strong result for the lower bound of the restricted $L^2$ norm is known \cite{prepare}:
\[
||\phi||_{L^2(\gamma)}\gg_\epsilon \tau^{-\epsilon}.
\]
Therefore, by Lemma \ref{key}, we have the second part of Theorem \ref{theorem}. However, as shown in \cite{pa}, a much weaker lower bound
\begin{equation}\label{iq}
||\phi||_{L^2(\gamma)}\gg e^{-B\tau},
\end{equation}
for some constant $B$, yields the same conclusion. In \cite{pa}, they call such a curve on a surface satisfying this lower bound ``good." The main result of this paper can be rephrased as:
\begin{itemize}
\item[(1)] any closed horocycle on a hyperbolic surface of finite volume is ``good,'' and
\item[(2)] any geodesic circle on a compact hyperbolic surface is ``good.''
\end{itemize}
We obtain such lower bounds by using purely analytic and geometric techniques. The main observation for our proof is Lemma \ref{cola}. From this, to deduce Theorem \ref{theorem}, we note that  it suffices to find a large Fourier coefficient ($b_\phi(n)$ as in Equation (\ref{exp1}) and $a_\phi(n)=b_{\phi,1}(n)$ in Equation (\ref{exp2})) in the range $c_1\tau < |n| <c_2 \tau$ for certain constants $c_1$ and $c_2$.

When $\gamma$ is a geodesic circle and $Y$ is compact, one can study the partial summation $S(X)=\sum_{|n|<X} |b_\phi (n)|^2$ by integrating $|\phi|^2$ on a region bounded by two concentric geodesic circles that contains multiple copies of $Y$. By quantifying the lower bound and the upper bound of $S(X)$, we get
\[
S(c_2 \tau) - S(c_1\tau) \gg 1
\]
for appropriate constants $c_1$ and $c_2$ (Lemma \ref{lemmageod}).

However, when $Y$ is non-compact and $\gamma$ is either a geodesic circle or a closed horocycle, the same method only yields $|N(\phi) \cap \gamma| \ll \tau \log^2 \tau$. Here the discrepancy between $\ll \tau \log^2 \tau$ and the desired bound $\ll \tau$ is closely related to the hypothetical phenomenon, known as ``escape of mass'' of the eigenfunctions. Hence we have to use a different method for non-compact surfaces.

When $\gamma$ is a closed horocycle, we use the theory of Eisenstein series and their analytic continuation. By integrating $|\phi|^2$ against an Eisenstein series and using the functional equation for the Eisenstein series, we can study the Fourier coefficients $a_\phi(n)$ by using the approximate functional equation. After choosing an appropriate test function in the approximate functional equation, we deduce Lemma \ref{lemmahoroc}.

\section{Special functions}
We begin our discussion with the uniform approximation formulas for some special functions. Let us define
\[
\tilde{K}_{i\tau}(x) = \frac{1}{|\Gamma(1/2+i\tau)|}K_{i\tau}(x) \text{ (}x>0\text{)}
\]
and
\[
C_\tau(m,x) = \frac{|\Gamma(1/2+m+i\tau)|}{|\Gamma(1/2+i\tau)|} P_{-1/2+i\tau}^{-m}(x) \text{ (} x>1\text{)}
\]
for a nonnegative integer $m$. Here $K_\tau(z)$ is the modified Bessel function of the second kind and $P_\tau^m(z)$ is the associated Legendre function of the first kind.
\begin{lem}[\cite{er}]\label{bessel}
We have
\begin{itemize}
\item[\rm{1}.] When $x>\tau>0$, put $\cos \theta = \tau/x$ where $\theta \in (0,\pi/2)$.
\begin{align*}
\tilde{K}_{i\tau} (x) & = \frac{\sqrt{\pi/2}+O(1/x)}{(x^2 -\tau^2)^{1/4}} e^{-(\sin \theta - \theta \cos \theta)x}
\end{align*}
\item[\rm{2}.] When $\tau>x>0$, put $\psi(\tau,x)=\tau\cosh^{-1}(\tau/x)-\sqrt{\tau^2-x^2}+ \pi/4$.
\begin{align*}
\tilde{K}_{i\tau} (x) & = \frac{\sqrt{2\pi}}{(\tau^2 -x^2)^{1/4}}\big{(}\sin(\psi(\tau,x))+O(\min\{\frac{1}{x},\frac{x^2}{\tau}\})\big{)}
\end{align*}
\item[\rm{3}.] When $|\tau-x|<\tau^{1/3}$, put $C_0 = 2^{-2/3} 3^{-1/6}\pi \Gamma(1/3)$.
\begin{align*}
\tilde{K}_{i\tau} (x) \approx C_0 x^{-1/3} + O(x^{-2/3})
\end{align*}
\end{itemize}
\end{lem}
\begin{lem}\label{legendre}
We have
\begin{itemize}
\item[\rm{1}.]
When $\beta:=\frac{\tau}{m} \ll 1$, $m>0$, and $x>1$ (\cite{dun}, p.325)
\[
C_\tau(m,x)= (\frac{2}{\pi})^{\frac{1}{2}}\big(\frac{\beta^2 - \eta}{x^2\beta^2-1-\beta^2}\big)^{\frac{1}{4}}\tilde{K}_{i\tau}(m\sqrt{\eta})\big(1+O(\frac{\sqrt{\eta}}{m})\big).
\]
Here the variable $\eta$ is given by:
\begin{align*}
\xi &= \frac{1}{x^2 -1}\\
\int_{\beta^2}^{\eta} \frac{(s-\beta^2 )^{\frac{1}{2}}}{2s} ds& = \int_{\beta^2}^{\xi} \frac{(s-\beta^2 )^{\frac{1}{2}}}{2s(1+s)^{\frac{1}{2}}} ds \text{ when } \xi > \beta^2\\
\int_{\beta^2}^{\eta} \frac{(\beta^2-s)^{\frac{1}{2}}}{2s} ds& = \int_{\beta^2}^{\xi} \frac{(\beta^2-s )^{\frac{1}{2}}}{2s(1+s)^{\frac{1}{2}}} ds \text{ when } \xi < \beta^2
\end{align*}
\item[\rm{2}.]
When $\alpha:= \frac{m}{\tau} \ll 1$, $m\geq 0$, and $x>1$ (\cite{dun}, p.314)
\[
C_\tau(m,x)= \big(\frac{\zeta-\alpha^2}{x^2-\alpha^2-1}\big)^{\frac{1}{4}} J_m(\tau \sqrt{\zeta})(1+O(\frac{1}{\tau})).
\]
Here the variable $\zeta$ is given by:
\begin{align*}
f(\alpha,s)&=\frac{1+\alpha^2-s^2}{(s^2-1)^2}\\
x_\alpha &= \sqrt{1+\alpha^2}\\
\int_{\alpha^2}^\zeta \frac{(s - \alpha^2)^{1/2}}{2s}ds &= \int_{x_\alpha}^x (-f(\alpha,s))^{1/2} ds \text{ when }x>x_\alpha\\
\int_{\alpha^2}^\zeta \frac{(\alpha^2-s)^{1/2}}{2s}ds &= \int_{x_\alpha}^x {f(\alpha,s)}^{1/2} ds \text{ when }x<x_\alpha
\end{align*}
\end{itemize}
\end{lem}
From these formulas we observe:
\begin{lem}\label{cola}
For any fixed $x>0$, we can find a constant $c_1>0$ with the following property:

For any constant $c_2>c_1$, there exists a constant $A>0$ depending on $c_2$ such that
\[
C_\tau(m,\cosh x) > e^{-Am}
\]
for any $c_1 \tau < m < c_2 \tau$, for all sufficiently large $\tau>0$.

The same statement is true for $\tilde{K}_{i\tau}(mx)$.
\end{lem}
\section{Analysis on hyperbolic surfaces}
\subsection{Upper half plane}
In order to study eigenfunctions on the hyperbolic surface $Y$, we identify $Y$ with a quotient space $\Gamma\backslash \mathbb{H}$ for some discrete subgroup $\Gamma \subset SL(2,\mathbb{R})$. Here $\mathbb{H}=\{z=x+iy|y>0\}$ is the upper half plane with the line element $ds^2=y^{-2}(dx^2+dy^2)$. Let $\triangle = \triangle_{\mathbb{H}}$ be the Laplace-Beltrami operator on $\mathbb{H}$ and let $dA=y^{-2}dxdy$ be the volume element.

By a function on $\Gamma \backslash \mathbb{H}$, we mean a function on $\mathbb{H}$ which is invariant under $\Gamma$. Let $\mathscr{F}$ be a fundamental domain for the action of $\Gamma$ on $\mathbb{H}$. For $1\leq p<\infty$, we let $L^p(\Gamma \backslash \mathbb{H})$ to be the space of measurable functions on $\Gamma \backslash \mathbb{H}$ such that
\[
||f||_{L^p(\Gamma \backslash \mathbb{H})}^p := \int_{\mathscr{F}}|f(z)|^p dA <\infty.
\]

We recall two theorems regarding eigenfunctions on the upper half plane.
\begin{lem}\label{geoexp}
Fix a point $z_0 \in \mathbb{H}$ and consider geodesic polar coordinate $(r,\theta)$ centered at $z_0$. For any bounded function $\phi$ on $\mathbb{H}$ that satisfies
\[
-\triangle \phi (z) = (1/4+\tau^2) \phi (z),
\]
there exist constants $c_\phi(m)$ such that
\[
\phi(r,\theta) = \sum_{m\in \mathbb{Z}} c_\phi(m) C_\tau(|m|, \cosh r) e^{im\theta}.
\]
\end{lem}
\begin{proof}
\cite{au}, p.142.
\end{proof}
\begin{lem}\label{expansion}
Assume $\phi$ is a function on $\mathbb{H}$ satisfying the conditions
\begin{align*}
-\triangle \phi (z) &= (1/4+\tau^2) \phi (z),\\
\phi(z)&=\phi(z+a),
\end{align*}
for some constant $a >0$. Assume further that there exist constants $C>0$ and $A>0$ such that
\[
|\phi(x+iy)| < Cy^A
\]
for $y>1$. Then there exist constants $\alpha_\phi$, $\beta_\phi$, and $c_\phi(m)$ such that
\[
\phi (x+iy) = \alpha_\phi y^{1/2+i\tau}+\beta_\phi y^{1/2-i\tau}+ \sum_{m \neq 0} c_\phi(m)  \sqrt{y}\tilde{K}_{i\tau}(2 \pi |m|y) e^{2 \pi i m x / a}.
\]
\end{lem}
\begin{proof}
\cite{au}, p.136.
\end{proof}
Note that in both lemmas, it is the growth condition of the eigenfunction which determines the choice of the special functions appearing in the equations. In general, without the growth condition, the associated Legendre function of the second kind and the modified Bessel function of the first kind will appear.

\subsection{Compact surface}
Let $Y = \Gamma \backslash \mathbb{H}$ be a compact hyperbolic surface. Fix a point $z_0\in \mathbb{H}$ and let $(r,\theta)$ be the geodesic polar coordinate centered at $z_0$. Note that any geodesic circle centered at $z_0$ is given by $\gamma=\{(r,\theta)|r=r_0\}$ for some $r_0$. Fix a fundamental domain $\mathscr{F}$ containing $z_0$ which has diameter equal to the diameter of $Y$.
\begin{lem}\label{lemma1}
Let $d$ be the diameter of $Y$ and let $f\in L^2(Y)$ be an $L^2(Y)$ normalized function. For any $a,b \in \mathbb{R}^+$ with $2d < b-a$,
\[
\int_a^b\int_0^{2\pi} |f(r,\theta)|^2 \sinh r d \theta d r > \frac{2\pi}{\mathrm{Area}(Y)} (\cosh (b-d) - \cosh (a+d))
\]
\end{lem}
\begin{proof}
Let $S$ be the set of $\gamma \in \Gamma$ such that $\gamma\mathscr{F} \subset \{(r,\theta)|a<r<b\}$. Observing that $\{(r,\theta)|a+d<r<b-d\} \subset S\mathscr{F} $,
\begin{align*}
\int_a^b\int_0^{2\pi} |f(r,\theta)|^2 \sinh r d \theta d r &\geq \int_{S\mathscr{F}}|f(r,\theta)|^2 dA\\
&=\frac{1}{\mathrm{Area}(Y)}\int_{S\mathscr{F}} dA\\
&> \frac{1}{\mathrm{Area}(Y)}\int_{a+d}^{b-d}\int_0^{2\pi} \sinh r d\theta dr\\
&=\frac{2\pi}{\mathrm{Area}(Y)} (\cosh (b-d) - \cosh (a+d)).
\end{align*}
\end{proof}
Let $\phi$ be an eigenfunction for $\triangle$ on $\Gamma \backslash \mathbb{H}$ with eigenvalue $-(1/4+\tau^2)$. We assume that $\phi$ is $L^2(Y)$ normalized. Since $Y$ is compact, $\phi$ is a bounded function on $\mathbb{H}$. Hence we can find constants $b_\phi(m)$ such that
\begin{equation}\label{exp1}
\phi(r,\theta) = \sum_{m\in \mathbb{Z}} b_\phi(m) C_\tau(|m|, \cosh r) e^{im\theta}.
\end{equation}
\begin{lem}\label{lemmageo}
We have
\[
\sum_{|m|<X} |b_\phi(m)|^2 \ll X+\tau.
\]
\end{lem}
\begin{proof}
From the asymptotics of $C_\tau(m,x)$, one proves
\[
\int_1^R |C_\tau(|m|,x)|^2 dx \gg \frac{1}{\tau}
\]
for $|m| < \tau \sqrt{R^2-1}$ and $1\ll R$. Therefore,
\begin{align*}
B(R):&= \int_0^{\cosh^{-1}R} \int_0^{2\pi} |\phi(r,\theta)|^2 \sinh r d\theta dr\\
&=\sum_{m \in \mathbb{Z}} |b_\phi(m)|^2\int_1^R |C_\tau(|m|,x)|^2 dx\\
&\gg \frac{1}{\tau}\sum_{|m| < \tau \sqrt{R^2-1}}|b_\phi(m)|^2
\end{align*}
for $R>\tau$. By using the same argument in Lemma \ref{lemma1}, we prove that $B(R) \ll R$. Therefore
\[
\sum_{|m| \leq \tau \sqrt{R^2-1}}|b_\phi(m)|^2 \ll \tau R,
\]
and substituting $\tau \sqrt{R^2-1}=X$ gives
\begin{align*}
\sum_{|m| \leq X}|b_\phi(m)|^2 &\ll \sqrt{X^2+\tau^2}\\
&<X+\tau .
\end{align*}
\end{proof}

\subsection{Noncompact surface and Eisenstein series}
In this section we state the results regarding noncompact hyperbolic surfaces and Eisenstein series. For a detailed treatment of the subject, see \cite{co}, \cite{ku}, and \cite{se}.

Assume that $Y=\Gamma \backslash \mathbb{H}$ is a noncompact hyperbolic surface which has finite volume. Let $\{\kappa_1 , \cdots, \kappa_h\}$ be a complete set of inequivalent cusps of $\Gamma$. Note that if $\kappa$ is a cusp of $\Gamma$, then there exists a $\sigma \in SL(2,\mathbb{R})$ such that $\sigma \kappa = \infty$ and $\sigma^{-1} \Gamma_\kappa \sigma = \Gamma_0$. Here $\Gamma_\kappa \subset \Gamma$ is the stabilizer subgroup of $\kappa$ and $\Gamma_0 = \{\left( \begin{array}{cc} 1 & n \\ 0 & 1 \end{array} \right)| n \in \mathbb{Z}\}$. This is equivalent to saying that $\sigma^{-1} \Gamma \sigma$ is reduced at $\infty$. We fix our choice of such $\sigma_j$ for each $\kappa_j$. Note that any closed horocycle on $Y$ is given by
\[
\{\sigma_j z|\mathrm{Im}(z)=y_0\}
\]
for some fixed $y_0>0$ and for some $1\leq j\leq h$.

Let $\phi(z)$ be a $L^2(Y)$ normalized real eigenfunction corresponding to the eigenvalue $-(\frac{1}{4}+\tau^2)$ with $\tau \geq 0$. Say the Fourier expansion of $\phi$ at $\kappa_j$ is given by
\begin{equation}\label{exp2}
\phi(\sigma_j z) = \sum_{n\neq 0} b_{\phi,j}(n) \sqrt{y} \tilde{K}_{i\tau}(2\pi|n|y) e^{2\pi i n x}
\end{equation}
for some constants $b_{\phi,j}(n)$ for each $j=1,\cdots, h$. Note that the constant terms in Lemma \ref{expansion} do not appear since $\phi \in L^2(Y)$ and $\tau \geq 0$. Also, because $\phi(z)$ is assumed to be real, $b_{\phi,j}(n) = \overline{b_{\phi,j}(-n)}$ for all $n \neq 0$.
\begin{lem}\label{lemma horo}
For any $j\in \{1,\cdots, h\}$,
\[
\sum_{|n|<X} |b_{\phi,j}(n)|^2 \ll X+\tau.
\]
\end{lem}
\begin{proof}
This is a slight generalization of \cite{si}, p.316, using the discreteness of $\Gamma$ and the Bruhat decomposition of $\Gamma$.
\end{proof}

Now let $E_j(z,s)$ be the Eisenstein series corresponding to $\kappa_j$.
\begin{thm}\label{mero}
Let $\Phi(s)=\left(\varphi_{ij}(s)\right)$ be the constant term matrix of $E_i(z,s)$. Then each $E_i(z,s)$ has a meromorphic continuation to the whole $s$-plane and the column vector $\mathscr{E}(z,s) = {^t}(E_1,\cdots,E_h)$ satisfies the functional equation
\[
\mathscr{E}(z,s) = \Phi(s) \mathscr{E}(z,1-s).
\]
\end{thm}
We remark that each $\varphi_{ij}(s)$ can be represented as the ratio of two holomorphic functions of finite order (\cite{eisen}, p.89 and \cite{se}), and that on $\mathrm{Re}(s)\geq \frac{1}{2}$, the functions $\varphi_{ij}(s)$ and $E_i$ have only finite number of poles lying in $(\frac{1}{2},1]$. Therefore we may fix a set of holomorphic functions $f(s),f_1(s), \cdots, f_h(s)$ of finite order such that
\begin{equation}\label{equation}
f(s)E_1(z,s) = \sum_{j=1}^h f_j(1-s)E_j(z,1-s)
\end{equation}
with $f(s)E_1(z,s)$ having a simple pole only at $s=1$ and $f_j(s)E_j(z,s)$ holomorphic in the region $\mathrm{Re}(s) \geq \frac{1}{2}$.

Following the standard Rankin-Selberg method, we obtain
\begin{lem}\label{R-S}
Define $L_j(s,\phi) = \sum_{n=1}^\infty \frac{|b_{\phi,j}(n)|^2}{n^s}$ whenever the series is absolutely convergent. There exist constants $e_j \neq 0$ depending only on the choice of $\sigma_j$ such that
\[
\int_Y E_j (z,s) |\phi (z)|^2 dA = e_j \cosh \pi \tau F(s,\tau) L_j(s,\phi)
\]
where
\[
F(s,\tau)=\pi^{-s} \frac{\Gamma^2(\frac{s}{2})\Gamma(\frac{s}{2}-i\tau)\Gamma(\frac{s}{2}+i\tau)}{\Gamma(s)}.
\]
\end{lem}
From Lemma \ref{lemma horo}, the defining series of $L_j (s,\phi)$ is absolutely convergent and $L_j(s,\phi) \ll |\zeta(s)|+\tau$ for $\mathrm{Re}(s) >1$. Note that $L_j(s,\phi)$ has a meromorphic continuation to the whole complex plane by Theorem \ref{mero}.

For the rest of the section, we study the growth of $F(\sigma+it,\tau)L_j(\sigma+it,\phi)$ at infinity for fixed $\sigma \geq \frac{1}{2}$ and ``fixed $\phi$." The only aim of this discussion is to take care of the technical convergence issues for the integrations in the proof of Lemma \ref{theorem ap}.

Using the Maass-Selberg relation, we have
\begin{lem}\label{lemma blah}
Fix a fundamental domain $Y_j$ for the action of $\sigma_j^{-1} \Gamma \sigma_j$ on $\mathbb{H}$ so that it is contained in a strip of width $1$ and contains a strip $Y_j(X)$ determined by $\mathrm{Im}(z)>X$ for sufficiently large $X>0$. Define $E_j^X$ by
\begin{align*}
E_j^X(\sigma_j z,s) :=
& \left\{
\begin{array}{cl}E_j(\sigma_j z,s) - y^s &\hspace{10mm} \text{ if }\mathrm{Im}(z)>X,\\
E_j(\sigma_j z,s) & \hspace{10mm}\text{ otherwise }.
\end{array} \right.
\end{align*}
for $z$ in the fundamental domain. Then we have the inner product formula
\[
(E_j^X(\sigma_j z,s),E_j^X(\sigma_j z,\bar{s'}))=\frac{X^{s+s'-1}}{s+s'-1}+ \varphi_{jj}(s') \frac{X^{s-s'}}{s-s'}-\varphi_{jj}(s)\frac{X^{-s+s'}}{s-s'}
\]
where $\mathrm{Re}(s)>1$, $\mathrm{Re}(s')>1$, and $s \neq s'$.
\end{lem}
Because of the location of the poles of the Eisenstein series, this inner product formula is valid for all $s$ and $s'$ in the region $\{s\in \mathbb{C}|\mathrm{Re}(s)>1/2, \text{ }s \notin (1/2,1]\}$ (p.32 \cite{ku}). Since $E_j^X(\cdot,s) \in L^2(Y_j)$ (p.23 \cite{co}) for any $s$ in such a region, taking the limit $s'\to s$ on the right hand side yields:
\[
(E_j^X(\sigma_j z,s),E_j^X(\sigma_j z,\bar{s}))=\frac{X^{2s-1}}{2s-1}+2\log X \varphi_{jj}(s) + \varphi_{jj}'(s).
\]
Let us denote the characteristic function of $Y_j(X)$ by $\chi_X(z)$. Then
\begin{align*}
\int_{Y_j} E_j (\sigma_j z,s) |\phi (\sigma_j z)|^2 dA &= \int_{Y_j} (E_j^X (\sigma_j z,s)+\chi_X(z)y^s) |\phi (\sigma_j z)|^2 dA\\
&\ll ||E_j^X(\cdot,s)||_{L^2(Y_j)} ||\phi||_{L^4(Y_j)}^2+\int_{Y_j(X)}y^s|\phi (\sigma_j z)|^2 dA.
\end{align*}
Note that $\phi$ decays exponentially in $y$ and $\varphi_{jj}(\sigma+it)$ has moderate growth at infinity for any fixed $\sigma > 1/2$(p.28 \cite{ku}). Then by Lemma \ref{lemma blah}, we see that $F(\sigma+it,\tau)L_j(\sigma+it,\phi)$ grows no faster than $e^{|t|}$ at infinity for any fixed $\sigma >1/2$. For $\sigma=1/2$, by \cite{eigen} we have
\[
\int_{Y_j} E_j (\sigma_j z,\frac{1}{2}+it) |\phi (\sigma_j z)|^2 dz \ll_\phi |t|^B e^{-\frac{\pi|t|}{2}}
\]
for some $B \in \mathbb{R}$. To summarize:
\begin{lem}\label{summary}
For any fixed $\sigma\geq 1/2$ and for $|t|>1$,
\[
F(\sigma+it,\tau)L_j(\sigma+it,\phi) \ll_\phi e^{|t|}
\]
for $j=1,\cdots, h$.
\end{lem}

\section{Complexification}
Let $D_\epsilon$ be an annulus on the complex plane given by $\{q~|~e^{-\epsilon} <|q|< e^\epsilon \}$. We first recall a lemma from \cite{pa}.
\begin{lem}\label{1}
Let $h$ be a holomorphic function defined on $D_\epsilon$. Normalize $h$ so that
\[
\int_0^1 |h(e^{2\pi ix})|^2 \mathrm{d}x=1.
\]
Then for any $\epsilon >0$ there exists a constant $C(\epsilon)>0$ such that the number of zeros of $h$ in $D_{\epsilon/2}$, which we denote by $n_\epsilon (h)$ is bounded by
\[
C(\epsilon)\max_{z \in D_\epsilon} \log |h(z)|.
\]
\end{lem}
Fix $\gamma$ to be either a geodesic circle or a closed horocycle. Say we are given the Fourier expansion for an eigenfunction on $\gamma$ so that
\[
\phi|_\gamma (x) = \sum_{n \in \mathbb{Z}} a_\phi (n) e^{i n x}.
\]
From the asymptotics of special functions and by Lemmas \ref{lemmageo} and \ref{lemma horo}, in either case we can find some constant $c>1$ depending only on $\gamma$ such that
\begin{align*}
|a_\phi (n)| \ll \tau &\text{ if } |n|<c\tau,\\
|a_\phi (n)| \ll e^{-|n|} &\text{ if } |n|\geq c\tau,
\end{align*}
where the implied constants depend only on $c$ and $\gamma$.

Define $\psi(q)$ by the relations $q=e^{i x}$ and $\psi(q)=\phi |_\gamma (x)$. If $q$ is in $D_{1/2}$, then we have
\begin{align*}
|\psi(q)| &\ll \sum_{|n|<c\tau} \tau e^{|n|/2} + \sum_{|n|>c\tau} e^{-|n|/2}\\
&\ll e^{c\tau}.
\end{align*}
Therefore $\psi(q)$ is holomorphic on $D_{1/2}$, and applying Lemma \ref{1} to $\psi(q)$, we obtain
\begin{lem}\label{key}
Let $Y$ be a hyperbolic surface which is compact (or has finite volume.) Let $\gamma$ be a fixed geodesic circle (or a fixed closed horocycle.) Then we have
\[
|N(\phi) \cap \gamma| \ll_{Y,\gamma} \log \frac{e^{c\tau}}{||\phi||_{L^2(\gamma)}}
\]
for some $c>1$ depending only on $Y$ and $\gamma$. Here $||f||_{L^2(\gamma)}$ is the restricted $L^2$ norm of $f$ on $\gamma$.
\end{lem}
Now, we are going to prove $||\phi||_{L^2(\gamma)}$ is not too small, from which Theorem \ref{theorem} will follow.

\section{Geodesic circle}
Assume $Y=\Gamma\backslash\mathbb{H}$ is compact and $\gamma$ is a fixed geodesic circle. Without loss of generality, we assume
\[
\gamma=\{(r,\theta)|r=r_0\}
\]
in the geodesic polar coordinate centered at some point $z_0 \in \mathbb{H}$. Let $b_\phi(m)$ be constants such that
\[
\phi(r,\theta) = \sum_{m\in \mathbb{Z}}b_\phi(m)C_\tau(|m|,\cosh r) e^{2\pi i m \theta}
\]
in the same coordinate system.

\begin{lem}\label{lemmageod}
For any $c_1>0$, there exists $c_2$ greater than $c_1$ such that for all sufficiently large $\tau$, there exists $m$ such that $c_1 \tau<|m|<c_2 \tau$, which satisfies $\frac{1}{\tau} \ll |b_\phi(m)|^2$.
\end{lem}
\begin{proof}
If $\alpha = \frac{m}{\tau} \leq c_1$, then by Lemma \ref{legendre},
\[
|C_\tau(m,x)|^2 \ll \big(\frac{\zeta-\alpha^2}{x^2-\alpha^2-1}\big)^{\frac{1}{2}} J_m(\tau \sqrt{\zeta})^2.
\]
From the asymptotics of Bessel function(\cite{er}), for any $x$ such that $\sqrt{\zeta}> 2 \alpha$,
\[
|C_\tau(m,x)|^2 \ll \big(\frac{\zeta-\alpha^2}{x^2-\alpha^2-1}\big)^{\frac{1}{2}} \frac{1}{\tau \sqrt{\zeta}}.
\]
Note that, because $\alpha$ is bounded by $c_1$, $\sqrt{\zeta}> 2 \alpha$ is ensured when $x$ is greater than some sufficiently large constant $C_1$ depending only on $c_1$. From the defining equation of $\zeta$, we observe that $\zeta < x^2 -1$ holds for any $\zeta>\alpha^2$. Since $\frac{\zeta-\alpha^2}{\zeta}$ is an increasing function in $\zeta$, we have
\begin{align*}
|C_\tau(m,x)|^2 &\ll \frac{1}{\tau \sqrt{x^2 - 1}} \\
&\ll \frac{1}{x\tau}.
\end{align*}
Now let $R>C_1$ and $a>1$ be constants which will be determined later. In Lemma \ref{legendre}, note that if $x$ is bounded from above by $aR$, then $\xi >\frac{1}{a^2R^2-1}$. Pick $c_2$ large enough so that
\[
\frac{1}{a^2R^2-1}>\frac{1}{c_2^2}>0
\]
and
\[
\int_{\frac{1}{c_2^2}}^{\frac{1}{a^2R^2-1}}\frac{(s-\frac{1}{c_2^2})^{\frac{1}{2}}}{2s(1+s)^{\frac{1}{2}}} ds > \frac{1}{c_2} \int_1^4\frac{(s-1)^{\frac{1}{2}}}{2s} ds.
\]
Then it follows from the defining equation of $\eta$ that $\eta > 4 \beta^2$ for any $ \beta \leq \frac{1}{c_2}$. We also deduce from the same equation that there exists some constant $\epsilon >0$ depending only on $a$ and $R$ such that $\eta>\epsilon$. Therefore, if $x < aR$ and $m \geq c_2 \tau$, then the asymptotics of $\tilde{K}_{i\tau}(x)$ imply
\[
|C_\tau(m,x)|^2 \ll e^{-Bm}
\]
for some constant $B>0$, depending only on $a$ and $R$.

Now set $|d_\phi (m)|^2 = |b_\phi(m)|^2 + |b_\phi(-m)|^2$ for $m>0$, and $d_\phi(0) = b_\phi(0)$. We have
\begin{align*}
B(a,R):&=\int_{\cosh^{-1}R}^{\cosh^{-1}aR} \int_0^{2\pi}|\phi(r,\theta)|^2 \sinh r d\theta dr \\
&= \sum_{m=0}^{\infty} |d_\phi (m)|^2 \int_{\cosh^{-1}R}^{\cosh^{-1}aR} |C_\tau(m, \cosh r)|^2 \sinh r dr\\
&=\sum_{m=0}^{\infty} |d_\phi (m)|^2 \int_R^{aR} |C_\tau(m, x)|^2 dx.
\end{align*}
We split this summation into three pieces and apply above estimations to obtain
\[
B(a,R) \ll \frac{\ln a}{\tau}\sum_{m \leq c_1 \tau} |d_\phi (m)|^2 + R(a-1)\sum_{c_1 \tau< m<c_2 \tau} |d_\phi (m)|^2 + \sum_{c_2 \tau \leq m} |d_\phi (m)|^2 e^{-Bm},
\]
where we used the trivial bound $C_\tau(m,x) \ll 1$ for the second summation. Applying Lemma \ref{lemmageo}, we see that
\[
B(a,R) \ll \ln a+ R(a-1)\sum_{c_1 \tau< m<c_2 \tau} |d_\phi (m)|^2 + \sum_{c_2 \tau \leq m} m e^{-Bm}.
\]
On the other hand, by Lemma \ref{lemma1}, we have
\begin{align*}
B(a,R)&>\frac{2\pi}{\mathrm{Area}(\mathscr{F})} (\cosh (\cosh^{-1}(aR)-d) - \cosh (\cosh^{-1}(R)+d))\\
&\gg aR.
\end{align*}
Therefore by choosing $a$ and $R$ large enough and then $\tau$ large, we find $b_\phi(m)$ with the desired property.
\end{proof}
Note that Parseval's theorem implies
\[
||\phi||_{L^2(\gamma)}^2 = \sum |b_\phi(m)|^2 |C_\tau(m,\cosh r_0)|^2.
\]
Hence by Lemma \ref{cola} and Lemma \ref{lemmageod}, we obtain
\[
||\phi||_{L^2(\gamma)} \gg_{Y,\gamma} e^{-B\tau}
\]
for some constant $B>0$. Therefore from Lemma \ref{key}, we deduce the first case of Theorem \ref{theorem}.
\section{Closed horocycle}
Let $Y=\Gamma \backslash \mathbb{H}$ and let $\gamma$ be a fixed horocycle. With notation as in Section 3.3, we assume without loss of generality that $\Gamma$ is reduced at $\infty$(that is, $\kappa_1=\infty$ and $\sigma_1=Id$) and that $\gamma$ is given by
\[
\gamma = \{x+iy|y=y_0\}
\]
for some $y_0>0$. Set $a_\phi (n)=b_{\phi,1}(n)$.
\subsection{Approximate functional equation}
First note that Theorem \ref{mero}, Lemma \ref{R-S}, and Equation (\ref{equation}) imply
\begin{equation}\label{equation2}
f(s)F(s,\tau)L_1(s,\phi) = \sum_{j=1}^h e_jf_j(1-s)F(1-s,\tau)L_j(1-s,\phi).
\end{equation}
By our choice of $f,f_1,\cdots,f_h$, on $\mathrm{Re}(s)\geq\frac{1}{2}$, $f(s)L_1(s,\phi)$ has a simple pole at $s=1$ and $f_j(s)L_j(s,\phi)$ is holomorphic. With this functional equation, we introduce a variation of the approximate functional equation(p.98 \cite{iw2}).
\begin{lem}\label{theorem ap}
For a sufficiently large integer $m$, define $\psi$ by
\[
\psi(X)= \frac{1}{2\pi i}\int_{(2)} f(s) F(s,\tau)s(1+\frac{s^2}{4\tau^2}) X^{-s} e^{-s^{4m}} ds
\]
Then for $X >0$, we have
\[
\sum_{n =1}^\infty |a_\phi (n)|^2 \psi(\frac{n}{X}) = (1+\frac{1}{4\tau^2})\frac{c}{\cosh \pi \tau}X +\frac{\tau^2}{\cosh \pi \tau}O(\frac{1}{X}),
\]
for some constant $c$ depending only on $f(s)$.
\end{lem}
\begin{proof}
For simplicity, set
\begin{align*}
\Lambda(s,\phi) &= f(s)F(s,\tau)L_1(s,\phi)\\
\Lambda_j(s,\phi)&=e_jf_j(s)F(s,\tau)L_j(s,\phi)\\
H_m(s,\tau)&=s(1+\frac{s^2}{4\tau^2})e^{-s^{4m}}
\end{align*}
By shifting the contour and using the functional equation (Equation (\ref{equation2})), we obtain
\begin{align*}
 \sum_{n =1}^\infty |a_\phi (n)|^2 &\psi(\frac{n}{X})=\frac{1}{2\pi i}\int_{(2)} \Lambda(s)H_m(s,\tau) X^s ds\\
=&(1+\frac{1}{4\tau^2})\frac{c}{\cosh \pi \tau}X +\frac{1}{2\pi i}\int_{(\frac{1}{2})}\Lambda(s) H_m(s,\tau)X^s ds\\
=&(1+\frac{1}{4\tau^2})\frac{c}{\cosh \pi \tau}X + \sum_{j=1}^h\frac{1}{2\pi i}\int_{(\frac{1}{2})}\Lambda_j (s)H_m(1-s,\tau)X^{1-s}ds\\
=&(1+\frac{1}{4\tau^2})\frac{c}{\cosh \pi \tau}X + \sum_{j=1}^h\frac{1}{2\pi i}\int_{(2)}\Lambda_j (s)H_m(1-s,\tau)X^{1-s}ds\\
=&(1+\frac{1}{4\tau^2})\frac{c}{\cosh \pi \tau}X+R(\phi)
\end{align*}
for some constant $c$ depending only on $f(s)$. From
\begin{align*}
L_j(2+it,\phi) &\ll \tau,\\
F(2+it,\tau) &\ll \frac{\tau}{\cosh \pi \tau} e^t,
\end{align*}
we get $R(\phi) = \frac{\tau^2}{\cosh \pi \tau}O(\frac{1}{X})$.

We complete the proof by observing that we can choose $m$ large enough to make all the summations and the integrations converge absolutely, since every holomorphic function that appears in the computation is of finite order (by our choice of $f,f_1,\cdots,f_h$ and Lemma \ref{summary}).
\end{proof}
\begin{lem}\label{11}
For $X>1$, we have
\[
\psi(X) \ll \frac{1}{\cosh \pi \tau}\frac{1}{X\log^2 X}
\]
and for $X>0$, we have
\[
\psi(X) \ll \frac{X}{\cosh \pi \tau}.
\]
\end{lem}
\begin{proof}
The first assertion follows after shifting the contour to $\mathrm{Re}(s)=1$ and using integration by parts. For the second assertion, note that $ F(s,\tau)s(1+\frac{s^2}{4\tau^2})$ has no pole for $\mathrm{Re}(s)>-2$. Now shift the contour to $\mathrm{Re}(s)=-1$.

In both cases we used
\[
\frac{\partial^k}{\partial t^k} \left(F(\sigma+it,\tau)\right) \ll_{k,\sigma} \frac{\tau^{\max\{\sigma-1 , 0\}}}{\cosh \pi \tau} e^t \text{ (}\sigma>-2, \text{ }k\geq 0 \text{)}
\]
which follows from Stirling's formula.
\end{proof}
\subsection{Lower bound for $||\phi||_{L^2(\gamma)}$}
\begin{lem}\label{lemmahoroc}
For any $c_1>0$, there exists $c_2$ greater than $c_1$ such that for all sufficiently large $\tau$, we can find $a_\phi (n)$ in range $c_1 \tau< n <c_2 \tau$ which satisfies $1 \ll |a_\phi (n)|^2$.
\end{lem}
\begin{proof}
From Lemma \ref{theorem ap}, we have
\[
\sum_{n =1}^\infty |a_\phi (n)|^2 \psi(\frac{n}{X}) = (1+\frac{1}{4\tau^2})\frac{c}{2 \cosh \pi \tau}X +\frac{\tau^2}{\cosh \pi \tau}O(\frac{1}{X}).
\]
Let $c_2>1$ and $\beta$ be the constants to be determined such that $0<c_1 < \beta < c_2$. Put $X= \beta \tau$. Split the summation into the three ranges
\[
\sum_{n =1}^\infty |a_\phi (n)|^2 \psi(\frac{n}{X}) = \sum_{n < c_1 \tau} +\sum_{c_1 \tau \leq n \leq c_2 \tau} +\sum_{ c_2 \tau< n}.
\]
Call these sums $I_1,I_2$, and $I_3$ respectively. For $I_1$, we use Lemma \ref{lemma horo} and Lemma \ref{11} to obtain
\[
|I_1| \ll \frac{c_1(c_1+1)}{\beta}\frac{\tau}{\cosh \pi \tau} .
\]
For $I_2$, we use Lemma \ref{11} to obtain
\[
|I_2| \ll \frac{1}{\cosh \pi \tau}\frac{c_2}{c_1}\sum_{c_1 \tau \leq n \leq c_2 \tau} |a_\phi (n)|^2.
\]
For $I_3$, we use Lemma \ref{lemma horo}, Lemma \ref{11} and the summation by parts to obtain
\begin{align*}
|I_3| &\ll \sum_{c_2 \tau < n} \frac{|a_\phi (n)|^2}{\cosh \pi \tau} \big( \frac{n}{X}\log^2 \frac{n}{X} \big)^{-1}\\
&\ll \frac{\tau}{\cosh \pi \tau} \frac{\beta}{\log c_2 - \log \beta}.
\end{align*}
These estimates imply that there exists a fixed constant $C>0$ such that
\[
\sum_{c_1 \tau \leq n \leq c_2 \tau} |a_\phi (n)|^2 > \frac{c_1}{c_2}(\frac{\beta}{C} - \frac{1}{\beta} - \frac{c_1(c_1+1)}{\beta} - \frac{\beta}{\log c_2  - \log \beta})\tau
\]
uniformly in $c_1$, $c_2$, $\beta$, and $\tau$. Therefore for any given $c_1>0$, we can choose $\beta$ and then $c_2$ so that
\[
\sum_{c_1 \tau \leq n \leq c_2 \tau} |a_\phi (n)|^2 \gg_{c_1} \tau
\]
uniformly in $\tau$.
\end{proof}
By Parseval's theorem, we have
\[
||\phi||_{L^2(\gamma)}^2= \sum_{n \neq 0} |a_\phi(n)|^2 y_0 \tilde{K}_{i\tau}(2\pi |n| y_0)^2.
\]
Hence by Lemma \ref{cola} and Lemma \ref{lemmahoroc}, we find a constant $B>0$ such that
\[
||\phi||_{L^2(\gamma)} \gg_{Y,\gamma} e^{-B\tau}.
\]
Therefore, from Lemma \ref{key}, we deduce the second case of Theorem \ref{theorem}.
\section{Appendix}
\subsection{Equidistribution of the geodesic circles}
Let $Y$ be a compact hyperbolic surface and fix a point $z_0 \in Y$. In this section we prove that the geodesic circle centered at $z_0$ becomes equidistributed on $Y$ when we increase the radius of the circle. Accordingly, for any fixed continuous function $f$ on $Y$, we prove
\begin{equation}\label{equi}
\lim_{r \to \infty} \frac{1}{2\pi} \int_0^{2\pi} f(\theta,r) d\theta =\int_Y f dA
\end{equation}
where $(r,\theta)$ is the geodesic polar coordinate centered at $z_0$. Since $f$ is continuous and $Y$ is compact, $f$ is a bounded function, hence in $L^2(Y)$. Therefore by the spectral expansion of $f$ in terms of the Laplacian eigenfunctions on $Y$, it is sufficient to prove Equation \ref{equi} for the eigenfunctions.

Assume $f \in L^2(Y)$ satisfies
\[
\triangle f = s(1+s)f
\]
with $s(1+s)< 0$. Then by Lemma \ref{geoexp},
\[
\frac{1}{2\pi} \int_0^{2\pi} f(\theta,r) d\theta = c_fP_s^0(\cosh r)=c_fP_s(\cosh r)
\]
for some constant $c_f$, and by the orthogonality of the eigenfunctions,
\[
\int_Y f dA=0.
\]
Therefore Equation \ref{equi} is equivalent to
\[
\lim_{r \to \infty} P_s(\cosh r)=0.
\]
When $s = -1/2+i\tau$ for some $\tau>0$, this follows from Lemma \ref{legendre} and when $-1<s<0$, we use the integral representation (p.172, \cite{le}):
\[
P_s(\cosh r) = \frac{1}{\pi}\int_0^\pi (\cosh r + \sinh r \cos \theta)^s d\theta.
\]

\bigskip
\footnotesize
\noindent\textit{Acknowledgments.}
We thank Peter Sarnak for introducing the problem and many helpful discussions. We also thank Sung-Jin Oh, Arul Shankar, and Jacob Tsimerman for the careful reading of the manuscript and helpful comments.

\end{document}